\documentclass{birkjour}
 \usepackage{mathabx}
 \usepackage{hyperref}
 \usepackage[dvipsnames]{xcolor}
 \newtheorem{thm}{Theorem}[section]
 \newtheorem{cor}[thm]{Corollary}
 \newtheorem{lem}[thm]{Lemma}
 \newtheorem{prop}[thm]{Proposition}
 \theoremstyle{definition}
 \newtheorem{defn}[thm]{Definition}
 \theoremstyle{remark}
 \newtheorem{rem}[thm]{Remark}
 
 \numberwithin{equation}{section}

\begin{document}

%-------------------------------------------------------------------------
% editorial commands: to be inserted by the editorial office
%
%\firstpage{1} \volume{228} \Copyrightyear{2004} \DOI{003-0001}
%\seriesextraline{This is the Concrete Title of this Book\br H.E. R and S.T.C. W, Eds.}
%\firstpage{1}

\title[The Rhaly Operators on Köthe Spaces]
 {{The Rhaly Operators on Köthe Spaces}}

%----------Author 1
\author{Nazl\i\;Do\u{g}an}
\address{Fatih Sultan Mehmet Vak{\i}f University, 
    34445 Istanbul, Turkey}
\email{ndogan@fsm.edu.tr}

%\thanks{This work was completed with the support of our
%----------Author 2
%\author{A Second Author}
%\address{The address of\br
%the second author\br
%sitting somewhere\br
%in the world}
%\email{dont@know.who.knows}


%\date{January 1, 2004}

\begin{abstract} 
We introduce and study the Rhaly operator on Köthe spaces, with a primary focus on understanding its well-definedness, continuity, and compactness. We especially examine operators acting on power series spaces of both infinite and finite type. In the sequel, we provide integral representations for the Rhaly operator on the space of entire functions $H(\mathbb{C})$ and the space of holomorphic functions on the unit disc $H(\mathbb{D})$.
We also investigate the topologizability and power boundedness of the Rhaly operators, which leads to findings about their mean ergodicity, uniform mean ergodicity, and Ces\`{a}ro boundedness.
\end{abstract}

%%% ----------------------------------------------------------------------
\maketitle
%%% ----------------------------------------------------------------------

%%%%%%%%%%%%%%%%%%%%%%%%%%%%%%%%
%%%%%%%%%%%%%%%%%%%%%%%%%%%%%%%%%%%%%%%%%%%%%%%%%%%%%%%%%%%%%%%%%%%%%%%%%
\section{Introduction}
The Ces\`{a}ro operator $\mathcal{C}$ on the Hilbert space $\ell^2$, 
defined by
$$(\mathcal{C}a)(n)=\frac{1}{n+1}\sum^{n}_{k=0}a_k, \qquad n=0,1,2,\ldots,$$
has been a central object of study in operator theory since the foundational 
work of Brown, Halmos, and Shields \cite{BHS}, who established boundedness 
with $\|\mathcal{C}\|_{\ell^2}=2$ and identified the spectrum 
$\sigma(\mathcal{C})=\{z:|z-1|\leq 1\}$. 
The operator $\mathcal{C}$ is represented, with respect to the standard 
orthonormal basis $(e_n)_{n\in\mathbb{N}}$, by the Ces\`{a}ro matrix 
\begin{equation}
\displaystyle \renewcommand{\arraystretch}{1.2} 
\setlength{\arraycolsep}{4pt}
\begin{bmatrix}
1 &0&0&0&\cdots \\
\frac{1}{2}&\frac{1}{2}&0&0&\cdots\\
\frac{1}{3}&\frac{1}{3}&\frac{1}{3}&0&\cdots\\
\frac{1}{4}&\frac{1}{4}&\frac{1}{4}&\frac{1}{4}&\cdots\\
\vdots&\vdots&\vdots&\vdots&\ddots
\end{bmatrix} 
\renewcommand{\arraystretch}{1.0}
\setlength{\arraycolsep}{4pt}.
\end{equation} The Ces\`{a}ro operator has been extensively studied on various spaces; 
see, e.g., \cite{ABR5,ABR4,ABR8,BURS,BHS,GP,Ross} and the references therein. Investigating the Ces\`{a}ro operator on Fr\'echet spaces is also a topic of interest; see \cite{ABR6, ABR7, E1, E2, E3}. In particular, Albanese, Bonet, 
and Ricker \cite{ABR7} investigated $\mathcal{C}$ on power series spaces 
of finite type $\Lambda_{1}(\alpha)$, establishing power boundedness and 
uniform mean ergodicity when $\Lambda_{1}(\alpha)$ is nuclear.

The (discrete) generalized Ces\`{a}ro operators $\mathcal{C}_{t}$, $t\in \left[0,1\right]$ were first studied by Rhaly in \cite{R1, R2}. The operator $\mathcal{C}_{t}$ acts on  $\ell^{2}$ by
$$(\mathcal{C}_{t}a)(n)=\frac{1}{n+1}\sum^{n}_{k=0} t^{n-k}a_{k}$$
for $n=0,1,2,\dots$.
With respect to the standard orthonormal basis $(e_n)_{n\in\mathbb{N}}$, 
$\mathcal{C}_t$ is represented by the generalized Ces\`{a}ro matrix
\begin{equation}
\renewcommand{\arraystretch}{1.4}
\setlength{\arraycolsep}{4pt}
\begin{bmatrix}
1 & 0 & 0 & 0 & \cdots \\
\frac{t}{2} & \frac{1}{2} & 0 & 0 & \cdots \\
\frac{t^2}{3} & \frac{t}{3} & \frac{1}{3} & 0 & \cdots \\
\frac{t^3}{4} & \frac{t^2}{4} & \frac{t}{4} & \frac{1}{4} & \cdots \\
\vdots & \vdots & \vdots & \vdots & \ddots
\end{bmatrix}
\renewcommand{\arraystretch}{1.0}
\setlength{\arraycolsep}{5pt}
\end{equation}
where $t\in[0,1]$. Rhaly \cite{R1} showed that $\mathcal{C}_t$ is bounded on $\ell^2$ for 
all $t\in[0,1]$, compact for $t\in[0,1)$, and determined its spectrum. 
Generalized Ces\`{a}ro operators have since been studied on various spaces; 
see \cite{ABR9,ABR8} and the references therein.

In \cite{R3}, Rhaly introduced the infinite matrices
\begin{equation}
\begin{bmatrix}
a_{0} &0&0&0&\cdots \\
a_{1}&a_{1}&0&0&\cdots\\
a_{2}&a_{2}&a_{2}&0&\cdots\\
a_{3}&a_{3}&a_{3}&a_{3}&\cdots\\
\vdots&\vdots&\vdots&\vdots&\ddots
\end{bmatrix}
\end{equation} 
now known as Rhaly matrices. The classical Ces\`{a}ro matrix is a special 
case, corresponding to $\displaystyle a_n=\frac{1}{n+1}$. Rhaly \cite{R3} proved that the operator $R_a$ associated with this matrix 
is bounded on $\ell^2$ if and only if the limit $L=\lim_{n\to\infty}(n+1)|a_n|$ 
exists, and compact when $L=0$.

While Rhaly operators have been studied on Hilbert spaces, 
their behavior on Fréchet spaces — and in particular on Köthe 
spaces — has not been systematically investigated. This is the gap the 
present paper addresses. Köthe spaces provide a natural and flexible 
framework that encompasses power series spaces, 
and spaces of analytic functions.

In this paper, we introduce and systematically investigate Rhaly operators on Köthe spaces. Our main contributions are organized as follows. We first establish necessary and sufficient conditions for the well-definedness, continuity, and compactness of $R_\theta:K(a_{n,k})\to K(b_{n,k})$ in Section~3. Applied to power series spaces, this characterization yields a complete picture. Specifically, Theorem \ref{Thm1} demonstrates that the operator $R_\theta:\Lambda_\infty(\alpha) \to\Lambda_\infty(\alpha)$ is well-defined, continuous, and compact for every symbol $\theta\in\Lambda_\infty(\alpha)$, whereas Theorem \ref{Thm3} establishes that compactness on $\Lambda_1(\alpha)$ is equivalent to $\theta$ belonging to the dual space $(\Lambda_1(\alpha))'$, provided that $\Lambda_1(\alpha)$ is nuclear. Motivated by Jasiczak's integral representation of Toeplitz operators \cite{J1,J2}, Section~4 is devoted to deriving integral representations of Rhaly operators on $H(\mathbb{C})$ and $H(\mathbb{D})$ by exploiting the canonical topological isomorphisms $H(\mathbb{C})\cong\Lambda_\infty(n)$ and $H(\mathbb{D})\cong\Lambda_1(n).$ Section~5 addresses topologizability, 
power boundedness, and ergodic properties. We prove that $R_\theta$ is always $m$-topologizable and power bounded on $\Lambda_\infty(\alpha)$ for any $\theta\in\Lambda_\infty(\alpha)$, and provide explicit criteria for the same properties on $\Lambda_1(\alpha).$ In conjunction with the theorems of Albanese, Bonet, and Ricker \cite{ABR} and Kalmes and Santacreu \cite{K1}, these results allow us to deduce the mean ergodicity and Ces\`{a}ro boundedness of Rhaly operators on power series spaces, and consequently on $H(\mathbb{C})$ and $H(\mathbb{D}).$ This work complements the author's previous studies on Toeplitz and Hankel operators within the framework of Köthe spaces \cite{N1,N2,N3}.

%%%%%%%%%%%%%%%%%%%%%%%%%%%%%%%%%%%%%%%%%%%%%%%%%%%%%%%%%%%%%%%%%%%%%%%%%
\section{Preliminaries}
This section provides the essential background for the paper. We refer the reader to \cite{MV} for more details. 

A Fréchet space is a complete Hausdorff locally convex space whose topology can be described by a countable system of seminorms $(\|\cdot\|_{k})_{k\in \mathbb{N}}$.

Köthe spaces hold a prominent place within the theory of Fréchet spaces due to their canonical basis structure. A matrix $(a_{n,k})_{k,n\in \mathbb{N}}$ of non-negative numbers is called a Köthe matrix if the following conditions hold:
\begin{itemize}
\item[1.] For each $n\in \mathbb{N}$ there exists a $k\in \mathbb{N}$ with $a_{n,k}>0$. 
\item[2.] $a_{n,k}\leq a_{n,k+1}$ for each $n,k\in \mathbb{N}$.
\end{itemize}
%{\color{red} Buradaki Köthe uzaylarının tanımı başka şekilde verilebilir, yani $\lamba^{\infty}$ eklenebilir.}
Given a Köthe matrix $(a_{n,k})_{n,k\in \mathbb{N}}$, the space
\[ K(a_{n,k})=\bigg\{ x=(x_{n})_{n\in \mathbb{N}}\; \bigg|\;\|x\|_{k}:=\sum^{\infty}_{n=1}|x_{n}|a_{n,k}< \infty \;\text{for all}\;k\in \mathbb{N} \bigg\}\]
is called a Köthe space. Every Köthe space is a Fr\'echet space given by the semi-norms
in its definition. The dual space of a Köthe space $K(a_{n,k})$ is
\[ (K(a_{n,k}))^{\prime}=\bigg\{ y=(y_{n})_{n\in \mathbb{N}}\; \bigg|\; \sup_{n\in \mathbb{N}}|y_{n}a_{n,k}^{-1}|<+\infty\; \text{for some}\; k\in \mathbb{N}\bigg\}\]
see Proposition 27.3 of \cite{MV} for more information.

Let $\alpha=\left(\alpha_{n}\right)_{n\in \mathbb{N}}$ be a non-negative increasing sequence with $\displaystyle \lim_{n\rightarrow \infty} \alpha_{n}=\infty$. The power series space of finite type associated with $\alpha$ is defined as
$$\Lambda_{1}\left(\alpha\right):=\left\{x=\left(x_{n}\right)_{n\in \mathbb{N}}: \;\left\|x\right\|_{k}:=\sum^{\infty}_{n=1}\left|x_{n}\right|e^{-{1\over k}\alpha_{n}}<\infty \textnormal{ for all } k\in \mathbb{N}\right\},$$
and the power series space of infinite type associated with $\alpha$ is defined as
$$\displaystyle \Lambda_{\infty}\left(\alpha\right):=\left\{x=\left(x_{n}\right)_{n\in \mathbb{N}}:\; \left\|x\right\|_{k}:=\sum^{\infty}_{n=1}\left|x_{n}\right|e^{k\alpha_{n}}<\infty \textnormal{ for all } k\in \mathbb{N}\right\}.$$
Power series spaces constitute an important class of Köthe spaces. They include the spaces of holomorphic functions on $\mathbb{C}^{d}$ and on the polydisc $\mathbb{D}^{d}$, where $\mathbb{D}$ denotes the unit disk in $\mathbb{C}$ and $d\in \mathbb{N}$:
\[\mathcal{H}(\mathbb{C}^{d})\cong \Lambda_{\infty}(n^{\frac{1}{d}}) \quad\quad\text{and} \quad\quad \mathcal{H}(\mathbb{D}^{d})\cong \Lambda_{1}(n^{\frac{1}{d}}).
\]
The nuclearity of a power series space of finite type $\Lambda_{1}(\alpha)$ and of infinite type $\Lambda_{\infty}(\alpha)$ is equivalent to conditions $\displaystyle\lim_{n\to \infty} \frac{\ln n}{\alpha_{n}}=0$ and $\displaystyle\sup_{n\in \mathbb{N}} \frac{\ln n}{\alpha_{n}}<+\infty$,
respectively.
%%%%%%%%%%%%%%%%%%%%%%%%%%%%%%%%%%%%%%%%%%%%%%%%%%%%%%%%%%%%%%%%%%%%%%%%
 
Let $E$ and $F$ be Fr\'echet spaces. A linear map $T:E\to F$ is continuous 
if for every $k\in \mathbb{N}$ there exists $p\in \mathbb{N}$ 
and $C_k>0$ such that 
\[\|Tx\|_k\leq C_k\|x\|_p\]
for all $x\in E$.  $\mathcal{L}(E,F)$ denotes the vector space of all continuous linear maps from $E$ to $F$. If $E=F$, we write $\mathcal{L}(E)$ instead of $\mathcal{L}(E,E)$ for brevity. The composition of an operator $T\in \mathcal{L}(E)$ with itself n-times is denoted by $T^{n}$. A linear map $T: E\to F$ is called compact
if $T(U)$ is precompact in $F$ where $U$ is a neighborhood of zero of E. 

Throughout this paper, $e_{n}$ stands for the sequence with 1 in the $n^{th}$ component and zeros elsewhere. In any Köthe space, the sequence $(e_{n})_{n\in \mathbb{N}}$ forms the canonical Schauder basis.

%%%%%%%%%%%%%%%%%%%%%%%%%%%%%%%%%%%%%%%%%%%%%%%%%%%%%%%%%%%%%%%%%%%%%%%%%%%%%%%%%%%%%%%%%%%%%%%%%%%%%%%%%%%%%%%%%%%%%%%%%%
To establish the continuity and compactness of operators defined between Köthe spaces, we apply  Lemma 2.1 of \cite{CR75}.

\begin{lem}[Lemma 2.1, \cite{CR75}]\label{Crone} 
	Let $K(a_{n,k})$ and $K(b_{n,k})$ be Köthe spaces.
	\begin{itemize}
		\item[a.] $T: K(a_{n,k}) \to K(b_{n,k})$ is a linear, continuous operator if and only if for each $k\in \mathbb{N}$, there exists $m\in \mathbb{N}$ such that
		\[ \sup_{n\in \mathbb{N}} {\frac{\|Te_{n}\|_{k}}{\|e_{n}\|_{m}}<\infty }.\]
		\item[b.] If $K(b_{n,k})$ is Montel, then $T: K(a_{n,k}) \to K(b_{n,k})$ is a compact operator if and only if there exists m such that for all k
		\[\sup_{n\in \mathbb{N}} {\frac{\|Te_{n}\|_{k}}{\|e_{n}\|_{m}}<\infty }.\]
	\end{itemize}
\end{lem}

The next proposition clarifies that the continuity condition on the basis elements is sufficient to ensure the well-definedness of the linear operator.

\begin{prop}\label{N1}\cite[Proposition 2.2]{N1}. Let $K(a_{n,k})$, $K(b_{n,k})$ be Köthe spaces and  $(a_{n})_{n\in \mathbb{N}}\subseteq K(b_{n,k})$  be a sequence. Let us define a linear map $T:K(a_{n,k})\to K(b_{n,k})$ by
	\[Te_{n}=a_{n} \hspace{0.5in}\text{and}\hspace{0.5in} Tx=\sum^{\infty}_{n=1}x_{n}Te_{n}\]
	for every $\displaystyle x=\sum^{\infty}_{n=1} x_{n}e_{n}$ and $n\in \mathbb{N}$. If the  continuity condition 
	\[\forall k\in \mathbb{N} \quad\quad \exists m\in \mathbb{N}\quad\quad\quad\quad \sup_{n\in \mathbb{N}}\frac{ \|Te_{n}\|_{k}}{\|e_{n}\|_{m}}<\infty \]
	holds, then $T$ is a well-defined and continuous operator.
\end{prop}

Let $T$ be a continuous linear operator on a Fr\'echet space $E$. $T^{[k]}$ denotes the k-th Ces\`{a}ro mean given by
 $$\frac{1}{k}\sum^{k}_{m=1} T^{m}$$
 for every $k\in \mathbb{N}$.
\begin{defn}\label{D1}Let $E$ be a Fr\'echet space. An operator $T\in \mathcal{L}(E)$ is called 
\begin{itemize}
\item \textbf{topologizable} if for every $p\in \mathbb{N}$ there exists $q\in \mathbb{N}$ such that for every $k\in \mathbb{N}$ there is $M_{k,p}>0$ satisfying
$$\|T^{k}(x)\|_{p}\leq M_{k,p}\|x\|_{q}$$
for every $x\in E$.
\item \textbf{m-topologizable} if for every $p\in \mathbb{N}$ there exist $q\in \mathbb{N}$ and $C_{p}\geq 1$ such that 
$$\|T^{k}(x)\|_{p}\leq C^{k}_{p}\|x\|_{q}$$
holds for every $k\in \mathbb{N}$ and $x\in E$.
\item \textbf{power bounded} if for every $p\in \mathbb{N}$ there exist $q\in \mathbb{N}$ and $C_{p}\geq 1$ such that
$$\|T^{k}(x)\|_{p}\leq C_{p}\|x\|_{q}$$
holds for every $k\in \mathbb{N}$ and $x\in E$.
\item \textbf{Ces\`{a}ro bounded }
if for every $p\in \mathbb{N}$ there exist $q\in \mathbb{N}$ and $C_{p}\geq 1$ such that 
$$\|T^{[k]}(x)\|_{p}\leq C_{p}\|x\|_{q}$$
holds for every $k\in \mathbb{N}$ and $x\in E$.
\end{itemize}
\end{defn}

It is straightforward that power bounded operators are m-topologizable and m-topologizable operators are topologizable.  In  \cite{B} Bonet proved that every bounded operator $T\in \mathcal{L}(E)$ on a locally convex space $E$ is m-topologizable \cite[Proposition 7]{B}. Therefore, every compact operator is m-topologizable.

By \cite[Lemma 2.5]{N3}, topologizability, $m$-topologizability, power 
boundedness, and Ces\`{a}ro boundedness of operators on Köthe spaces can be verified by examining the operator only on the canonical basis $(e_n)_{n\in \mathbb{N}}$. 
Throughout this paper, we apply this lemma to study these properties 
for Rhaly operators.

A continuous linear operator $T$ on a locally convex Hausdorff space $E$ is called \textbf{mean ergodic} if the limits
$$Px=\lim_{n\to \infty} \frac{1}{n}\sum^{n}_{m=1} T^{m}x, \hspace{0.15in} x\in E$$
exist in $E$. If the convergence is uniform on bounded subsets of E, then T is
called \textbf{uniformly mean ergodic}. A locally convex Hausdorff space E is defined as (uniformly) mean ergodic if every power bounded operator acting on it is (uniformly) mean ergodic. 

Albanese, Bonet, and Ricker showed in \cite[Proposition 2.13]{ABR} that 
a Köthe space is mean ergodic if and only if it is Montel. Since all power 
series spaces are Montel, every power bounded operator on a power series 
space is uniformly mean ergodic. Furthermore, Kalmes and Santacreu proved 
in \cite[Theorem 2.5]{K1} that on Montel spaces, mean ergodicity is 
equivalent to Ces\`{a}ro boundedness together with 
$\displaystyle \frac{T^n}{n}\to 0$ pointwise.

%%%%%%%%%%%%%%%%%%%%%%%%%%%%%%%%%%%%%%%%%%%%%%%%%%%%%%%%%%%%%%%%%%%%%%%%%%%
\section{Rhaly Operators}
Given a sequence $\theta=(\theta_{n})_{n\in \mathbb{N}_{0}}$, the associated Rhaly matrix is
\[
\begin{bmatrix}
\theta_{0} &0&0&0&\cdots \\
\theta_{1}&\theta_{1}&0&0&\cdots\\
\theta_{2}&\theta_{2}&\theta_{2}&0&\cdots\\
\theta_{3}&\theta_{3}&\theta_{3}&\theta_{3}&\cdots\\
\vdots&\vdots&\vdots&\vdots&\ddots
\end{bmatrix}.
\]
We introduce the operator $R_{\theta}: K(a_{n,k})\to K(b_{n,k})$ associated with the above matrix by defining its action on the canonical basis vectors as
\[R_{\theta}e_{n}=(0,\dots,0, \theta_{n-1},\theta_{n},\dots)=\sum^{\infty}_{j=n}\theta_{j-1}e_{j}\]
and extending it by linearity:
\begin{equation*}\label{E1} 
\quad R_{\theta}x=\sum^{\infty}_{n=1}x_{n} R_{\theta}e_{n}= \sum^{\infty}_{n=1}\left(\sum^{n}_{j=1} x_{j}\right)\theta_{n-1}e_{n}, \quad \text{for } x=\sum^{\infty}_{n=1}x_{n}e_{n}\in K(a_{n,k}).
\end{equation*}

Since the convergence of the series $\displaystyle \sum_{n=1}^\infty x_n R_\theta e_n$ in $K(b_{n,k})$ is not guaranteed for every $x \in K(a_{n,k})$, the operator $R_\theta$ is not well-defined in general. The aim of this section is to determine the conditions on $\theta$ and on the Köthe spaces that ensure $R_\theta$ is well-defined and continuous, with a particular emphasis on power series spaces.

If the operator $R_{\theta}:K(a_{n,k})\to K(b_{n,k})$ is well-defined, then $R_{\theta}(e_1)=\theta\in K(b_{n,k})$ must hold. This shows that assuming $\theta\in K(b_{n,k})$ is a natural condition. Moreover, the inequalities
\[
\|R_{\theta}e_{n}\|_{k}=\sum^{\infty}_{j=n}|\theta_{j-1}|b_{j,k} \leq \sum^{\infty}_{j=1}|\theta_{j-1}|b_{j,k}=\|\theta\|_{k}
\]
guarantee that $R_{\theta}e_{n}\in K(b_{n,k})$ for all $n\in \mathbb{N}$, provided that $\theta\in K(b_{n,k})$. Proposition \ref{N1} leads us to the following.

\begin{prop}\label{RP1} The operator $R_{\theta}:K(a_{n,k})\to K(b_{n,k})$ is well-defined and continuous if and only if $\theta\in K(b_{n,k})$ and the continuity condition 
\[\forall k\in \mathbb{N} \quad\quad \exists m\in \mathbb{N} \quad\quad\quad\quad \sup_{n\in \mathbb{N}}\frac{ \|R_{\theta}e_{n}\|_{k}}{\|e_{n}\|_{m}}<\infty \]
holds.
\end{prop}

Our initial focus is to characterize the conditions under which the operator $R_{\theta}$ defined between the Köthe spaces $K(a_{n,k})$ and $K(b_{n,k})$ is well-defined, continuous, and compact.

\begin{prop}\label{P2} Let $K(a_{n,k}), K(b_{n,k})
$ be  Köthe spaces and let $\theta\in K(b_{n,k})$. If $\displaystyle \inf_{n\in \mathbb{N}} a_{n,m_{0}}>0$ for some $m_{0}\in\mathbb{N}$, then the Rhaly operator $R_{\theta}:K(a_{n,k}) \to K(b_{n,k})$ is well-defined and continuous. If, in addition $K(b_{n,k})$ is Montel, then $R_{\theta}$ is compact.
\end{prop}
\begin{proof} Let us assume that $\theta\in K(b_{n,k})$ and that $\displaystyle \inf_{n\in \mathbb{N}} a_{n,m_{0}}\geq C>0$ for some $m_{0}\in\mathbb{N}$ and $C>0$ holds.  This implies that for every $k\in \mathbb{N}$ there exists a constant $D_{k}>0$ such that
\begin{equation*}
\|R_{\theta}e_{n}\|_{k}=\sum^{\infty}_{j=n}|\theta_{j-1}|b_{j,k}\leq \sum^{\infty}_{j=1}|\theta_{j-1}|b_{j,k} \leq \|\theta\|_{k}\leq \frac{1}{ C}\|\theta\|_{k}a_{n,m_{0}}\leq D_{k} a_{n,m_{0}}
\end{equation*}
for every $n\in \mathbb{N}$ and 
$$\sup_{n\in \mathbb{N}} \frac{\|R_{\theta}e_{n}\|_{k}}{\|e_{n}\|_{m_{0}}}<\infty.$$
Hence $R_{\theta}$ is well-defined and continuous by Proposition \ref{N1}. If $K(b_{n,k})$ is a Montel space, $R_{\theta}$ is also compact by Lemma \ref{Crone} as $m_{0}$ does not depend on $k$.
\end{proof}

%------------------------------------------------------------------------------------------------------------

As an immediate consequence of Proposition \ref{P2}, we obtain the following.

\begin{thm}\label{Thm1}  Let $K(b_{n,k})
$ be a Köthe space and let $\theta\in K(b_{n,k})$.  The Rhaly operator $R_{\theta}: \Lambda_{\infty}(\alpha) \to K(b_{n,k})$ is well-defined and continuous. Furthermore, if $K(b_{n,k})$ is Montel, then $R_{\theta}$ is also compact. In particular, the following results are obtained.
\begin{itemize}
\item[1.] If $\theta\in \Lambda_{\infty}(\beta)$, then the operator $R_{\theta}:\Lambda_{\infty}( \alpha)\to \Lambda_{\infty}(\beta)$ is well-defined, continuous and compact. 
\item[2.] If $\theta\in \Lambda_{1}(\beta)$, then the operator $R_{\theta}:\Lambda_{\infty}(\alpha)\to \Lambda_{1}(\beta)$ is well-defined, continuous and compact. 
\end{itemize}
\end{thm}

We now seek a weaker assumption on the matrix 
of $K(a_{n,k})$ in Proposition \ref{P2}. When $K(b_{n,k})$ is a power series space of infinite type $\Lambda_{\infty}(\beta)$, a stability condition on the sequence $\beta$ is required: we say that the sequence $\beta$ is weakly-stable if
$$
\sup_{n\in \mathbb{N}} {\beta_{n+1}\over \beta_{n}}<\infty.$$

\begin{prop}\label{R5}
Let $\beta=(\beta_{n})_{n\in \mathbb{N}}$ be a weakly stable sequence and let 
$\theta\in \Lambda_{\infty}(\beta)$. The Rhaly operator $R_{\theta}:K(a_{n,k})\to \Lambda_{\infty}(\beta)$ 
is 
\begin{itemize}
\item[1.] well-defined and continuous provided that
\begin{equation}\label{Eq3}
\forall k\in \mathbb{N} \quad \exists m\in \mathbb{N},\, C>0  \quad 
e^{-k\beta_{n}}\leq C\, a_{n,m} \qquad \forall n\in \mathbb{N}.
\end{equation}
\item[2.] compact provided that
\begin{equation}\label{REQ4}
\exists m\in \mathbb{N} \quad \forall k\in \mathbb{N} \quad \exists C>0  \qquad 
e^{-k\beta_{n}}\leq C\, a_{n,m} \quad \forall n\in \mathbb{N}.
\end{equation}
\end{itemize}
\end{prop}
\begin{proof} Since $\beta$ is a weakly-stable sequence, there exists an $M\in \mathbb{N}$, $M> 1$, satisfying
$$\beta_{n+1}\leq M\beta_{n}$$
for every $n\in \mathbb{N}$. Let us assume that $\theta\in \Lambda_{\infty}(\beta)$ and the condition (\ref{Eq3}) holds. Then, for every $k\in \mathbb{N}$ there exist $m\in \mathbb{N}$, $C_{1},C_{2}>0$ such that
\begin{align*}
\begin{split}
\|R_{\theta}e_{n}\|_{k}&=\sum^{\infty}_{j=n} |\theta_{j-1}|e^{k\beta_{j}} =\sum^{\infty}_{j=n} |\theta_{j-1}|e^{k\beta_{j}} e^{k\beta_{j+1}} e^{-k\beta_{j+1}}
\\ &\leq e^{-k\beta_{n+1}} \sum^{\infty}_{j=n} |\theta_{j-1}|e^{k(M+1)\beta_{j}} \leq e^{-k\beta_{n}} \sum^{\infty}_{j=1} |\theta_{j-1}|e^{k(M+1)\beta_{j}} \\
&\leq C_{1} a_{n,m}\|\theta\|_{(M+1)k}  
\leq C_{2}\|e_{n}\|_{m}.
\end{split}
\end{align*}
holds for every $n\in \mathbb{N}$ and
$$\label{T1}
\sup_{n\in \mathbb{N}}\frac{\|R_{\theta}e_{n}\|_{k}}{\|e_{n}\|_{m}}< \infty.
$$
Hence, $R_{\theta}$ is well-defined and continuous by  Proposition \ref{RP1}.

The compactness characterization in (2) follows by the same argument, 
with the index $m$ chosen independently of $k$, which ensures compactness 
by Lemma \ref{Crone}.
\end{proof}

Combining parts (1) and (2) of Proposition \ref{R5}, we conclude that the operator  $R_{\theta}:\Lambda_{1}(\alpha)\to \Lambda_{\infty}(\beta)$  is well-defined, continuous, and compact, provided certain conditions on $\alpha$ and $\beta$ are met.

\begin{thm}
Let $\alpha=(\alpha_n)_{n\in\mathbb{N}}$ and $\beta=(\beta_{n})_{n\in \mathbb{N}}$ 
be sequences with $\beta$ weakly stable, and let $\theta\in \Lambda_{\infty}(\beta)$. 
If $\displaystyle\sup_{n\in\mathbb{N}}\frac{\alpha_n}{\beta_n}<\infty$, then 
the Rhaly operator $R_{\theta}:\Lambda_{1}(\alpha)\to \Lambda_{\infty}(\beta)$ 
is well-defined, continuous, and compact.
\end{thm}
\begin{proof} Since $\displaystyle \sup_{n\in\mathbb{N}}\frac{\alpha_n}{\beta_n}<\infty$, there exists 
$A\in \mathbb{N}$ such that $\displaystyle \frac{\alpha_n}{\beta_n}\leq A$, that is, 
$\alpha_{n}\leq A\beta_{n}$ for all $n\in \mathbb{N}$.
Then for all $m,k\in \mathbb{N}$ we have
$$ \frac{1}{mA}\alpha_{n}\leq \frac{1}{A}\alpha_{n}\leq \beta_{n}\leq k\beta_{n}$$
and therefore 
$$-k\beta_{n}\leq -\frac{1}{mA}\alpha_{n} \qquad \text{and} \qquad e^{-k\beta_{n}}\leq e^{-\frac{1}{mA}\alpha_{n}}$$
for all $n\in \mathbb{N}$.
Setting $\tilde{m}=mA\in\mathbb{N}$, we obtain
$$e^{-k\beta_n}\leq e^{-\frac{1}{\tilde{m}}\alpha_{n}}$$
for all $n\in\mathbb{N}$. Since $\tilde{m}$ can be chosen 
independently of $k$, both conditions \eqref{Eq3} and \eqref{REQ4} 
are satisfied. By Proposition \ref{R5}, $R_{\theta}:\Lambda_{1}(\alpha)
\to \Lambda_{\infty}(\beta)$ is well-defined, continuous, and compact.
\end{proof}

We now turn to the case where the range space is a power series space of 
finite type $\Lambda_1(\beta)$ and establish the conditions for the Rhaly operator $R_\theta: K(a_{n,k}) \to \Lambda_1(\beta)$ 
to be well-defined, continuous, and compact.

\begin{prop}
\label{R2b}
Let $\theta\in \Lambda_{1}(\beta)$.  The Rhaly operator $R_{\theta}:K(a_{n,k})\to \Lambda_{1}(\beta)$ is
\begin{itemize}
\item[1.]
 well-defined and continuous provided that 
\begin{equation}\label{REQ3a}
\forall k\in \mathbb{N} \quad \exists m\in \mathbb{N},\, C>0  \qquad 
e^{-\frac{1}{k}\beta_{n}}\leq C\, a_{n,m} \qquad \forall n\in \mathbb{N}.
\end{equation}
\item[2.] compact provided that
\begin{equation}\label{REQ4a}
\exists m\in \mathbb{N} \quad \forall k\in \mathbb{N} \quad \exists C>0  \qquad 
e^{-\frac{1}{k}\beta_{n}}\leq C\, a_{n,m} \quad \forall n\in \mathbb{N}.
\end{equation}
\end{itemize}
\end{prop}

\begin{proof} Assume that $\theta\in \Lambda_{1}(\beta)$ and the condition (\ref{REQ3a}) holds. Then for every $k\in \mathbb{N}$, there exist $m\in \mathbb{N}$ and $C_{1}, C_{2}>0$ such that
\begin{align*}
\begin{split}
\|R_{\theta}e_{n}\|_{k}&=\sum^{\infty}_{j=n} |\theta_{j-1}|e^{-\frac{1}{k}\beta_{j}} =\sum^{\infty}_{j=n} |\theta_{j-1}|e^{-\frac{1}{2k}\beta_{j}}e^{-\frac{1}{2k}\beta_{j}} \leq  e^{-\frac{1}{2k}\beta_{n}}\sum^{\infty}_{j=n} |\theta_{j-1}|e^{-\frac{1}{2k}\beta_{j}} 
\\
&  \leq C_{1}a_{n,m} \sum^{\infty}_{j=n} |\theta_{j-1}|e^{-\frac{1}{2k}\beta_{j}}=C_{1}\|\theta\|_{2k}a_{n,m}\leq C_{2}a_{n,m}
\end{split}
\end{align*}
for every $n\in \mathbb{N}$, and therefore $$ \sup_{n\in \mathbb{N}}\frac{\|R_{\theta}e_{n}\|_{k}}{\|e_{n}\|_{m}}< \infty.$$
Hence $R:K(a_{n,k})\to \Lambda_{1}(\beta)$  is well defined and continuous by Proposition \ref{RP1}.

For compactness, assume that condition \eqref{REQ4a} holds. Then $m$ can 
be chosen independently of $k$, and the above estimates yield
$$\sup_{n\in \mathbb{N}}\frac{\|R_{\theta}e_{n}\|_{k}}{\|e_{n}\|_{m}}< \infty$$
with $m$ independent of $k$. Hence $R_{\theta}$ is compact by 
Lemma \ref{Crone}.
\end{proof}

%%%%%%%%%%%%%%%%%%%%%%%%%%%%%%%%%%%%%%%%%%%%%%%%%%%%%%%%%%%%%%%%%%%%%%%%%%%%%%%%

By Proposition \ref{R2b}, it follows that the operator $R_{\theta}:\Lambda_{1}(\alpha) \to \Lambda_{1}(\beta)$ is well-defined and continuous, given that specific conditions on $\alpha$ and $\beta$ hold.

\begin{thm}\label{Thm2} Let $\alpha=(\alpha_n)_{n\in\mathbb{N}}$ and $\beta=(\beta_{n})_{n\in \mathbb{N}}$ 
be sequences, and let $\theta\in \Lambda_{\infty}(\beta)$. 
If $\displaystyle\sup_{n\in\mathbb{N}}\frac{\alpha_n}{\beta_n}<\infty$, then 
the Rhaly operator $R_{\theta}:\Lambda_{1}(\alpha)\to \Lambda_{1}(\beta)$ is well-defined and continuous.
\end{thm}
\begin{proof}  Since $\displaystyle \sup_{n\in\mathbb{N}}\frac{\alpha_n}{\beta_n}<\infty$, there exists 
$A\in \mathbb{N}$ such that 
$\alpha_{n}\leq A\beta_{n}$ for all $n\in \mathbb{N}$. Then we have 
$$ \frac{1}{kA}\alpha_{n}\leq \frac{1}{k}\beta_{n}\hspace{0.25in} \text{and}\hspace{0.25in}-\frac{1}{k}\beta_{n}\leq -\frac{1}{kA}\alpha_{n}$$
for every $k,
n\in \mathbb{N}$.
Therefore, for all $k,m\in \mathbb{N}$ satisfying $m>kA$
$$e^{-\frac{1}{k}\beta_{n}}\leq e^{-\frac{1}{m}\alpha_{n}}$$
for every $n\in \mathbb{N}$. This says that the condition in (\ref{REQ3a}) is satisfied. $R_{\theta}:\Lambda_{1}(\alpha)\to \Lambda_{1}(\beta)$ is well-defined, continuous by Proposition \ref{R2b}.
\end{proof}

We now turn to the compactness of the Rhaly operator on a power series 
space of finite type. In this case, a complete characterization is possible.

\begin{thm}\label{Thm3} Let $\Lambda_{1}(\alpha)$ be a nuclear power series space of finite type  and $\theta\in \Lambda_{1}(\alpha)$. The Rhaly operator $R_{\theta}:\Lambda_{1}(\alpha)\to \Lambda_{1}(\alpha)$ is compact if and only if $\theta\in (\Lambda_{1}(\alpha))^{\prime}$.
\end{thm}
\begin{proof} Let us assume that the Rhaly operator $R_{\theta}:\Lambda_{1}(\alpha)\to \Lambda_{1}(\alpha)$ is compact, that is, there exists a  $m_{0}\in \mathbb{N}$ so that for every $p\in \mathbb{N}$, there is a $D_{p}>0$ satisfying
$$\sup_{n\in \mathbb{N}}\frac{\|R_{\theta}e_{n}\|_{p}}{\|e_{n}\|_{m_{0}}}\leq D_{p}$$
and this gives us that
$$|\theta_{n-1}|e^{\left(\frac{1}{m_0} -\frac{1}{p}\right)\alpha_{n}}\leq \frac{\|R_{\theta}e_{n}\|_{p}}{\|e_{n}\|_{m_{0}}}\leq D_{p}$$
for every $n\in \mathbb{N}$. We choose $p$ as $2m_0$. Then there exists a $D_{m_0}$ satisfying
$$|\theta_{n-1}|
\leq D_{m_0} e^{-\frac{1}{2m_0}\alpha_{n}}$$
for every $n\in \mathbb{N}$. This says that $\theta\in (\Lambda_{1}(\alpha))^{\prime}$. For the converse, let us assume that $\theta\in (\Lambda_{1}(\alpha))^{\prime}$. There exist $r_{0}\in \mathbb{N}$ and  $D_{r_{0}}>0$ so that
$$|\theta_{n-1}|\leq D_{r_{0}}e^{-\frac{1}{r_0}\alpha_{n}}$$
for every $n\in \mathbb{N}$. Then for every $p,n\in \mathbb{N}$, we have
\begin{equation*}
\begin{split}
\frac{\|R_{\theta}e_{n}\|_{p}}{\|e_{n}\|_{r_{0}}}&=e^{\frac{1}{r_{0}}\alpha_{n}} \sup_{j\geq n} |\theta_{j-1}| e^{-\frac{1}{p}\alpha_{j}} \leq  \sup_{j\geq n} |\theta_{j-1}| e^{(\frac{1}{r_0} -\frac{1}{p})\alpha_{j}} 
\\ & \leq D_{r_{0}} \sup_{j\geq n} e^{-\frac{1}{r_0} \alpha_{j}} e^{(\frac{1}{r_0} -\frac{1}{p})\alpha_{j}} \leq D_{r_{0}}.
\end{split}
\end{equation*}
This says that $R_{\theta}$ is compact.
\end{proof}

\begin{rem}
As shown by Albanese, Bonet, and Ricker in \cite[Proposition 2.5]{ABR7}, 
the Ces\`{a}ro operator on $\Lambda_{1}(\alpha)$ is not compact. When 
$\Lambda_{1}(\alpha)$ is nuclear, this is consistent with 
Theorem \ref{Thm3}, since the sequence 
$\left(\frac{1}{n+1}\right)_{n\in\mathbb{N}_{0}}$ does not belong to 
$(\Lambda_{1}(\alpha))^{\prime}$.
\end{rem}

%%%%%%%%%%%%%%%%%%%%%%%%%%%%%%%%%%%%%%%%%%%%%%%%%%%%%%%%%%%%%%%%%%%%%%%%%%%%%%%%%%%%%%%%%%%%%%%%%%%%%%%%%%%%%%%%%%%%%%%%
\section{Integral Representation of Rhaly Operators} 

In this section, motivated by Jasiczak's integral representation of Toeplitz operators on 
$H(\mathbb{C})$ and $H(\mathbb{D})$ \cite{J1,J2}, we adopt a similar perspective for Rhaly operators on $H(\mathbb{C})$ and $H(\mathbb{D})$. The representation of Toeplitz operators via Cauchy-type integrals 
suggests that operators defined through structured infinite matrices can be realized as integral operators with appropriately chosen kernels or symbols. 
Following this approach, we establish integral representations of the Rhaly 
operator on $H(\mathbb{C})$ and $H(\mathbb{D})$.

Throughout this section, let $G$ denote either $\mathbb{C}$ or $\mathbb{D}$. The space $H(G)$, 
endowed with the topology of uniform convergence on compact subsets, 
is isomorphic to the power series space
$$H(\mathbb{C})\cong\Lambda_{\infty}(n) \qquad \text{and} \qquad 
H(\mathbb{D})\cong\Lambda_{1}(n),$$
via the map $f\mapsto (a_n)_{n\in\mathbb{N}_0}$, where $a_n$ denotes 
the $n$-th Taylor coefficient of $f$ at the origin.

For every $g\in H(G)$, one defines the operator $R_{g}:H(G) \to H(G)$ such that for every $f\in H(G)$ and $z\in G$, put
\begin{equation}
(R_{g}f)(z)=\frac{1}{2\pi i}\int_{|w|=r_{0}} \frac{f(w)}{w(1-w)}g\left(\frac{z}{w}\right)dw
\end{equation}
where $r_{0}$ is chosen so that $|z|<r_{0}$ and $r_{0}\neq 1$. Since $f,g\in H(G)$ and the contour $|w|=r_{0}$ avoids the singularities 
$w=0$ and $w=1$, the function $\displaystyle w\to \frac{f(w)}{w(1-w)}g\left(\frac{z}{w}\right)$ is holomorphic on $|w|=r_{0}$. By Cauchy's 
theorem, the value of the integral is independent of the choice of $r_{0}$, 
and hence $R_g$ is well-defined.

\begin{thm}\label{TII}
Let $G$ denote either $\mathbb{C}$ or $\mathbb{D}$. Let $R: H(G)\to H(G)$ be a continuous linear operator. 
Then the following are equivalent:
\begin{itemize}
\item[(i)] The matrix of $R$ with respect to the Schauder basis 
$\{z^{n}\}_{n\in\mathbb{N}_{0}}$ is a Rhaly matrix, that is, there exists 
a sequence $\theta=(\theta_{n})_{n\in\mathbb{N}_{0}}$ such that
$$\hspace{1.6in}R(z^{m})=\sum^{\infty}_{j=m}\theta_{j-1}z^{j} \qquad \hspace{1.2in}m\in \mathbb{N}_{0}.$$
\item[(ii)] There exists a function $g\in H(G)$ such that $R=R_{g}$.
\end{itemize}
Moreover, in this case the sequence $\theta$ and the function $g$ are 
related by
$$\hspace{1.6in}\theta_{n}=\frac{1}{2\pi i} \int_{|w|=r} \frac{g(w)}{w^{n+1}}dw 
\qquad\hspace{1in} n\in \mathbb{N}_{0}$$
where r is a constant with $0<r<1$.
\end{thm}
\begin{proof} Let $\displaystyle f(z)=\sum^{\infty}_{n=0} a_{n}z^{n}$ and $\displaystyle g(z)=\sum^{\infty}_{n=0}b_{n}z^{n}$ be the elements of $H(G)$ We show that
$$(R_{g}f)(z)=\frac{1}{2\pi i}\int_{|w|=r_{0}} \frac{f(w)}{w(1-w)}g\left(\frac{z}{w}\right)dw=\sum^{\infty}_{n=0}\left(\sum^{n}_{k=0}a_{k}\right)b_{n}z^{n}$$
for every $z\in \mathbb{C}$. 

Let us define
$\displaystyle h(z)=\sum^{\infty}_{n=0}\left(\sum^{n}_{k=0}a_{k}\right)z^{n}$ for every $z\in \mathbb{D}$ and say $\displaystyle c_{n}=\sum^{n}_{k=0}a_{k}$ for every $n\in \mathbb{N}$. Then we have $\displaystyle h(z)=\sum^{\infty}_{n=0}c_n z^n$. By the following of equalities
\begin{equation*}
\begin{split}
h(z)&= \sum^{\infty}_{n=0} \left(\sum^{n}_{k=0}a_{k}\right)z^{n}=\sum^{\infty}_{k=0}a_{k}\left(\sum^{\infty}_{j=k} z^{j}\right) \\
&=\sum^{\infty}_{k=0} a_{k}z^{k}\left(\sum^{\infty}_{j=0}z^{j}\right)=\frac{1}{1-z}\sum^{\infty}_{k=0}a_{k}z^{k}=\frac{1}{1-z}f(z),
\end{split}
\end{equation*}
hence $\displaystyle h(z)=\frac{f(z)}{1-z}$ for every $z\in \mathbb{D}$. 

Fix a $z\in G$ and a 
$r_{0}$ satisfying $|z|<r_{0}$ and $r_{0}\neq 1$.  We have
\begin{equation*}
\begin{split}
&\frac{1}{2\pi i}\int_{|w|=r_{0}} \frac{f(w)}{w(1-w)}g\left(\frac{z}{w}\right)dw=\frac{1}{2\pi i}\int_{|w|=r_{0}} \frac{h(w)}{w}g\left(\frac{z}{w}\right)dw \\
&= \sum^{\infty}_{n=0}b_{n}z^{n} \frac{1}{2\pi i}\int_{|w|=r_{0}} \frac{h(w)}{w^{n+1}}dw=\sum^{\infty}_{n=0}c_nb_nz^n=\sum^{\infty}_{n=0}\left(\sum^{n}_{k=0}a_{k}\right)b_{n}z^n
\end{split}
\end{equation*}
for every $z\in \mathbb{C}$. Consequently,
$$(R_{g}z^{n})(\xi)=\sum^{\infty}_{m=n}b_{m}\xi^{m}$$
for every $n\in \mathbb{N}_{0}$. Hence, the associated matrix of $R_{g}$ is a Rhaly matrix and 
$$\theta_{n}=b_{n}=\frac{1}{2\pi i} \int_{|w|=1} \frac{g(w)}{w^{n+1}}dw$$
for every $n\in \mathbb{N}_{0}$.

On the other hand, if $R:H(G)\to H(G)$ is a continuous linear operator whose associated matrix is a Rhaly matrix, the function $\displaystyle g(z)=\sum^{\infty}_{n=0}\theta_n z^n$, formed by the coefficients of the sequence $\theta$, is an element of $H(G)$ since
\begin{equation*}
(R1)(z)=\sum^{\infty}_{n=0}\theta_{n}z^{n}=g(z)\in H(G). 
\end{equation*}
The above calculation likewise makes it evident that $R=R_{g}$.
\end{proof}

As $H(\mathbb{C})$ and $H(\mathbb{D})$ are isomorphic to $\Lambda_{\infty}(n)$ and $\Lambda_{1}(n)$, respectively, the following result is an immediate consequence of Theorem \ref{T1} and Theorem \ref{Thm3}:

\begin{cor} \begin{enumerate}
    \item The operator $R_{g}:H(\mathbb{C}) \to H(\mathbb{C})$ is compact for every $g\in H(\mathbb{C})$.
    \item The operator $R_{g}:H(\mathbb{D}) \to H(\mathbb{D})$ is compact if and only if $g$ belongs to the strong dual of $H(\mathbb{D})$.
\end{enumerate} 
\end{cor}

As an application, we identify the function $g\in H(\mathbb{D})$ that 
generates the Ces\`{a}ro operator via the integral representation established 
in Theorem \ref{TII}. Recall that the Ces\`{a}ro operator on $H(\mathbb{D})$ 
is defined by
$$(\mathcal{C}f)(z)=\frac{1}{z}\int^{z}_{0}\frac{f(\xi)}{1-\xi}d\xi, 
\qquad z\in \mathbb{D}.$$
It is well known that the matrix of $\mathcal{C}$ with respect to the Schauder basis $\{z^n\}_{n\in\mathbb{N}_0}$ is a Rhaly matrix. Indeed, one has
\[
\mathcal{C}(z^m)=\sum_{n=m}^{\infty}\frac{1}{n+1}z^n,
\qquad m\in\mathbb{N}_0,
\]
so that the associated Rhaly sequence is given by $\displaystyle \theta_n=\frac{1}{n+1}$.

By Theorem \ref{TII}, there exists a function $g\in H(\mathbb{D})$ such that $\mathcal{C}=R_g$, where the symbol $g$ is determined by the coefficients $\theta$. Hence,
\[
g(z)=\sum_{n=0}^{\infty}\frac{1}{n+1}z^n.
\]
Using the classical identity
\[
-\log(1-z)=\sum_{n=1}^{\infty}\frac{z^n}{n}, \qquad |z|<1,
\]
we obtain
\[
g(z)=
\begin{cases}
-\dfrac{\log(1-z)}{z}, & z\neq 0,\\[6pt]
1, & z=0,
\end{cases}
\]
where the singularity at $z=0$ is removable. In particular, $g\in H(\mathbb{D})$. Therefore, the Ces\`{a}ro operator corresponds to a Rhaly operator with symbol $g(z)=-\log(1-z)/z$.
%%%%%%%%%%%%%%%%%%%%%%%%%%%%%%%%%%%%%%%%%%%%%%%%%%%%%%%%%%%%%%%%%%%%%%%%%%%%%%%%%%%%%%%%%%%%%%%%%%%%%%%%%%%%%%%%%%%%%%%%%%%%%%%

\section{Topologizability and Power Boundedness of Rhaly Operators}

We establish conditions for the topologizability and power boundedness of 
the Rhaly operator on $\Lambda_{1}(\alpha)$, $\Lambda_{\infty}(\alpha)$, 
and certain Köthe spaces $K(a_{n,k})$. These results are then applied to 
determine the corresponding properties of Rhaly operators on $H(\mathbb{C})$ 
and $H(\mathbb{D})$.

In Proposition \ref{P2}, we showed that the operator $R_{\theta}:K(a_{n,k})\to K(b_{n,k})$ is compact if $\displaystyle \inf_{n\in \mathbb{N}} a_{n,m_{0}}>0$ for some $m_{0}\in \mathbb{N}$ and $K(a_{n,k})$ is Montel. As all compact operators are m-topologizable, the operator $R_{\theta}: K(a_{n,k})\to K(a_{n,k})$ is m-topologizable.  In the following proposition, we show that the operator $R_{\theta}: K(a_{n,k})\to K(a_{n,k})$ is m-topologizable without assuming that $K(a_{n,k})$ is Montel.

\begin{prop}\label{PBG} Let $K(a_{n,k})$ be a Köthe space and let $\theta\in K(a_{n,k})$. If  \\$\displaystyle \inf_{n\in \mathbb{N}}a_{n,m_{0}}>0$ for some $m_{0}\in \mathbb{N}$, then the Rhaly operator $R_{\theta}:K(a_{n,k})\to K(a_{n,k})$ is m-topologizable.
\end{prop}
\begin{proof} Let us assume that $\theta\in K(a_{n,k})$ and  $\displaystyle A=\inf_{n\in \mathbb{N}}  a_{n,m_{0}}>0$ holds for some $m_{0}\in\mathbb{N}$.  For every $p\in \mathbb{N}$, we define $\displaystyle D_{p}:=\frac{\|\theta\|_{p}}{A}\geq 0$. Then
\begin{equation*}
\|R_{\theta}e_{n}\|_{p}=\sum^{\infty}_{j=n}|\theta_{j-1}|a_{j,p}\leq \sum^{\infty}_{j=1}|\theta_{j-1}|a_{j,p} \leq \|\theta\|_{p}\leq\frac{1}{ A}\|\theta\|_{p}a_{n,m_{0}}= D_{p} \|e_{n}\|_{m_0}
\end{equation*}
for every $n\in \mathbb{N}$, and therefore
\begin{equation*}
\|R_{\theta}x\|_{p}\leq \sum^{\infty}_{n=1}|x_{n}|\|R_{\theta}e_{n}\|_{p}\leq D_{p}\sum^{\infty}_{n=1}|x_{n}|\|e_{n}\|_{m_{0}}\leq D_{p}\|x\|_{m_{0}}
\end{equation*}
for every $x\in K(a_{n,k})$. Moreover,
\begin{equation*}
\|R^{2}_{\theta}x\|_{p}\leq D_{p} \|R_{\theta}x\|_{m_{0}}\leq D_{p} D_{m_{0}} \|x\|_{m_{0}}   
\end{equation*}
for every $x \in \mathbb{N}$ and 
\begin{equation}\label{RGCPB}
\|R^{k}_{\theta}e_{n}\|_{p}\leq D_{p} (D_{m_{0}})^{k-1}\|e_{n}\|_{m_{0}} \leq (\max\{D_{p}, D_{m_{0}}\})^{k} \|e_{n}\|_{m_0}    
\end{equation}
for every $k,n\in \mathbb{N}$. By Lemma 2.5 of \cite{N3},  $R_{\theta}: K(a_{n,k})\to K(a_{n,k})$ is m-topologizable.
\end{proof}

As an immediate consequence of Proposition \ref{PBG}, we obtain the following.

\begin{cor}\label{IMT} The Rhaly operator $R_{\theta}:\Lambda_{\infty}(\alpha)\to \Lambda_{\infty}(\alpha)$ is m-topologizable for every $\theta\in \Lambda_{\infty}(\alpha)$.
\end{cor}

The following proposition provides a sufficient condition for the power 
boundedness of $R_{\theta}:K(a_{n,k})\to K(a_{n,k})$, analogous to 
Proposition \ref{PBG}.

\begin{prop}\label{PPG} Let $K(a_{n,k})$ be a Köthe space and let $\theta\in K(a_{n,k})$. If $\displaystyle A:=\inf_{n\in \mathbb{N}}a_{n,m_{0}}>0$ for some $m_{0}\in \mathbb{N}$ and $\displaystyle \|\theta\|_{m_{0}}\leq A$, then the Rhaly operator $R_{\theta}:K(a_{n,k})\to K(a_{n,k})$ is power bounded.
\end{prop}
\begin{proof} Let $A:=\inf_{n\in \mathbb{N}}a_{n,m_{0}}>0$ for some $m_0\in \mathbb{N}$ and  $\theta\in K(a_{n,k})$ with $ \|\theta\|_{m_{0}}\leq A$. Then the number $\displaystyle D_{m_{0}}=\frac{\|\theta\|_{m_{0}}}{A}$ is less than or equal to 1. By inequality (\ref{RGCPB}) and Lemma 2.5 of \cite{N3}, $R_{\theta}$ is power bounded. 
\end{proof}

Combining Proposition \ref{PPG} with  Proposition 2.13 of \cite{ABR} and Theorem 2.5 of \cite{K1}, 
we obtain the following.

\begin{cor}
 Let $K(a_{n,k})$ be a Montel Köthe space and let $\theta\in K(a_{n,k})$. If $\displaystyle A:=\inf_{n\in \mathbb{N}}a_{n,m_{0}}>0$ for some $m_{0}\in \mathbb{N}$ and $\displaystyle \|\theta\|_{m_{0}}\leq A$,  then the Rhaly operator $R_{\theta}:K(a_{n,k})\to K(a_{n,k})$ is 
\begin{itemize}
\item[i.] mean ergodic,
\item[ii.] uniformly mean ergodic, 
\item[iii.] Ces\`{a}ro bounded and $\displaystyle \lim_{n\to \infty} \frac{R^{n}_{\theta}}{n}=0,$ pointwise in $K(a_{n,k})$.
\end{itemize}
\end{cor}
%%%%%%%%%%%%%%%%%%%%%%%%%%%%%%%%%%%%%%%%%%%%%%%%%%%%%%%%%%%%%%%%%%%%%%%%%%%%%%%%%%%%%%%%%%%%%%%%%%%%%%%%

We now establish power boundedness for Rhaly operators on power series 
spaces of infinite type.

\begin{thm}\label{PBIT} The Rhaly operator $R_{\theta}:\Lambda_{\infty}(\alpha)\to \Lambda_{\infty}(\alpha)$ is power bounded for every $\theta\in \Lambda_{\infty}(\alpha)$.
\end{thm}
\begin{proof}
Let $\theta\in \Lambda_{\infty}(\alpha)$ and fix $p\in \mathbb{N}$.  Since 
$\|\theta\|_{p}<\infty$, there exists $q_p\in\mathbb{N}$ such that 
$\|\theta\|_p\leq e^{q_p\alpha_1}=\inf_{n\in\mathbb{N}}e^{q_p\alpha_n}$. 
By Proposition \ref{PPG}, $R_{\theta}:\Lambda_{\infty}(\alpha)\to
\Lambda_{\infty}(\alpha)$ is power bounded.
\end{proof}

Combining Theorem \ref{PBIT} with Proposition 2.13 of \cite{ABR} and Theorem 2.5 of \cite{K1}, 
we obtain the following.

 \begin{cor}\label{PBITC} For every $\theta\in \Lambda_{\infty}(\alpha)$,
 the Rhaly operator $R_{\theta}:\Lambda_{\infty}(\alpha)\to \Lambda_{\infty}(\alpha)$ is
\begin{itemize}
    \item[i.] mean ergodic,
    \item[ii.] uniformly mean ergodic,
    \item[iii.] Ces\`{a}ro bounded and $\displaystyle \lim_{n\to \infty} \frac{R_{\theta}^{\hspace{0.025in}n}}{n}=0,$ pointwise in $\Lambda_{\infty}(\alpha)$.
\end{itemize}
 \end{cor}

Since $H(\mathbb{C})$ and $\Lambda_{\infty}(n)$ are isomorphic, and all relevant properties 
are preserved under isomorphism, the following is obtained. 
For analogous results on other operators on power series spaces, 
see \cite{A5, K1}.

\begin{cor} For every $g\in H(\mathbb{C})$, the operator $R_{g}:H(\mathbb{C})\to H(\mathbb{C})$ is
\begin{itemize}
\item[i.] m-topologizable,
\item[ii.] power bounded,
\item[iii.] mean ergodic,
\item[iv.] uniformly mean ergodic,
\item[v.] Ces\`{a}ro bounded and $\displaystyle \lim_{n\to \infty} \frac{R^{n}_{g}}{n}=0,$ pointwise in $H(\mathbb{C})$.
\end{itemize}
\end{cor}
%%%%%%%%%%%%%%%%%%%%%%%%%%%%%%%%%%%%%%%%%%%%%%%%%%%%%%%%%%%%%%%%%%%%%%%%%%%%%%%%%%%%%%%%%%%%%%%%%%%%%%%%%%%%%%%%%%%%%%%%%%%%%%%%%%%5
Our next step is to explore the properties of the Rhaly operator $R_{\theta}:\Lambda_{1}(\alpha)\to \Lambda_{1}(\alpha)$. We start by computing its powers and their corresponding seminorms. 

Assume that $\theta\in \Lambda_{1}(\alpha)$ and $\Lambda_{1}(\alpha)$ is nuclear. Since for every $n\in \mathbb{N}$ and $x\in \Lambda_{1}(\alpha)$, $\displaystyle R_{\theta}e_{n}=\sum^{\infty}_{j=n} \theta_{j-1}e_{j}$ and $\displaystyle R_{\theta}x=\sum^{\infty}_{n=1} \left(\sum^{n}_{i=1} x_{i}\right)\theta_{n-1} e_{n}$,
$$R^{2}_{\theta}e_{n}=R_{\theta}(R_{\theta}e_{n})=\sum^{\infty}_{j=n} \theta_{j-1}\left(\sum^{j}_{i=n}\theta_{i-1}\right)e_{j} \hspace{0.15in} \text{and}\hspace{0.15in}(R^{2}_{\theta}e_{n})_m=\theta_{m-1}\sum^{m}_{j=n}\theta_{j-1}$$
for every $m\geq n$, $n\in \mathbb{N}$. Hence
$$\|R^{2}_{\theta}e_{n}\|_{p}\leq \sup_{m\geq n}|\theta_{m-1}|\left(\sum^{m}_{j=n}|\theta_{j-1}|\right)e^{-\frac{1}{p}\alpha_{m}}$$
for every $p,n\in \mathbb{N}$. In general, the coefficients of the power are written as
$$(R^{k}_{\theta}e_{n})_{m}=\sum_{n-1\leq j_{1}\leq j_{2}\leq \cdots\leq j_{k}=m-1}\theta_{j_{1}}\cdots \theta_{j_{k}}$$
for every $m\geq n$, and 0 for every $m<n$, $m,k,n\in \mathbb{N}$ . Hence, we obtain that
$$|(R^{k}_{\theta}e_{n})_{m}|\leq\left( \sum^{m}_{j=n}|\theta_{j-1}|\right)^{k}\hspace{0.15in}\text{and}\hspace{0.15in}
\|R^{k}_{\theta}e_{n}\|_{p}\leq \sup_{m\geq n}\left(\sum^{m}_{j=n}|\theta_{j-1}|\right)^{k}e^{-\frac{1}{p}\alpha_{m}}$$
for every $m\geq n$, $n,k,p\in \mathbb{N}$. Consequently, the inequalities 
\begin{equation} \label{Fesas}
\begin{split}
\|R^{k}_{\theta}e_{n}\|_{p}&\leq \sup_{m\geq n}\left(\sum^{m}_{j=n}|\theta_{j-1}|\right)^{k}e^{-\frac{1}{p}\alpha_{m}} \\
&=\sup_{m\geq n}\left(\sum^{m}_{j=n}|\theta_{j-1}|\right)^{k}e^{-\frac{1}{3p}\alpha_{m}}e^{-\frac{1}{3p}\alpha_{m}}e^{-\frac{1}{3p}\alpha_{m}} \\
&\leq e^{-\frac{1}{3p}\alpha_{n}} \sup_{m\geq n} \left(\sum^{m}_{j=n} |\theta_{j-1}|e^{-\frac{1}{3pk}\alpha_{j}}\right)^{k} e^{-\frac{1}{3p}\alpha_{m}} \\
& \leq e^{-\frac{1}{3p}\alpha_{n}}  \left(\sum^{\infty}_{j=n} |\theta_{j-1}|e^{-\frac{1}
{3pk}\alpha_{j}}\right)^{k} \sup_{m\geq n} e^{-\frac{1}{3p}\alpha_{m}} \\
&\leq e^{-\frac{1}{3p}\alpha_{n}} \left(\|\theta\|_{3pk}\right)^{k} =\|e_{n}\|_{3p} \left(\|\theta\|_{3pk}\right)^{k}
\end{split}
\end{equation}
are satisfied for every $n,k,p\in \mathbb{N}$.
%%%%%%%%%%%%%%%%%%%%%%%%%%%%%%%%%%%%%%%%%%%%%%%%%%%%%%%%%%%%%%%%%%%%%%%%%%%%%%%%%%%%%%

As an immediate consequence of inequality \eqref{Fesas} and 
Lemma 2.5 of \cite{N3}, we obtain the following.
\begin{prop} Let $\Lambda_{1}(\alpha)$ be a nuclear power series space of finite type and  $\theta\in \Lambda_{1}(\alpha)$. The Rhaly operator $R_{\theta}:\Lambda_{1}(\alpha)\to \Lambda_{1}(\alpha)$  is
\begin{itemize}
\item[i.] topologizable,
\item[ii.] m-topologizable if $\displaystyle \sup_{p\in \mathbb{N}}\|\theta\|_{p}=\sup_{n\in \mathbb{N}}|\theta_{n}|<\infty$,
\item[iii.] power bounded if $\displaystyle \sup_{p\in \mathbb{N}}\|\theta\|_{p}=\sup_{n\in \mathbb{N}}|\theta_{n}|\leq 1$, then $R_{\theta}$.
\end{itemize}
\end{prop}
%%%%%%%%%%%%%%%%%%%%%%%%%%%%%%%%%%%%%%%%%%%%%%%%%%%%%%%%%%%%%%%%%%%%%%%%%%%%%%%%%%%%%%%%%%%%%%%%%%%%
%\begin{thm}\label{FTP} Let $\Lambda_{1}(\alpha)$ be a nuclear space and  $\theta\in \Lambda_{1}(\alpha)$. The Rhaly operator $R_{\theta}:\Lambda_{1}(\alpha)\to \Lambda_{1}(\alpha)$ is power bounded if and only if $\displaystyle \sup_{p\in \mathbb{N}}\|\theta\|_{p}\leq 1$.
%\end{thm}
%\begin{proof}
%Let assume that $\displaystyle \sup_{p\in \mathbb{N}} \|\theta\|_{p}\leq 1$. (\ref{Fesas}) gives us that for every $p, k,n \in \mathbb{N}$, it yields that
%$$\|R^{k}_{\theta}e_{n}\|_{p}\leq \|e_n\|_{3p}.$$
%Then $R_{\theta}$ is power bounded.
%For the converse,  we assume that $R_{\theta}:\Lambda_{1}(\alpha)\to \Lambda_{1}(\alpha)$ is power bounded. This means that for every $p\in \mathbb{N}$ there exists a $q\in \mathbb{N}$ such that 
%$$\|R^{k}_{\theta}e_{n}\|_{p}\leq \|e_{n}\|_{q}$$
%holds for every $k,n\in \mathbb{N}$. By choosing $n=k=1$, we have
%$$\|\theta\|_{p}=\|R_{\theta}e_{1}\|_{p}\leq \|e_{1}\|_{q}=e^{-\frac{1}{q}\alpha_{1}}<1.$$
%This completes the proof.
%\end{proof}

As a direct consequence of Proposition 2.13 of \cite{ABR} and Theorem 2.5 of \cite{K1}, we establish the following result:

 \begin{cor}\label{PBITC1} Let $\Lambda_{1}(\alpha)$ be a nuclear power series space of finite type and $\theta\in \Lambda_{1}(\alpha)$. If $\displaystyle \sup_{p\in \mathbb{N}}\|\theta\|_{p}=\sup_{n\in \mathbb{N}}|\theta|_{n}\leq 1$, then the Rhaly operator $R_{\theta}:\Lambda_{1}(\alpha)\to \Lambda_{1}(\alpha)$ is
\begin{itemize}
    \item[1.] power bounded,
    \item[2.] mean ergodic,
    \item[3.] uniformly mean ergodic,
    \item[4.] Ces\`{a}ro bounded and $\displaystyle \lim_{n\to \infty} \frac{R_{\theta}^{\hspace{0.025in}n}}{n}=0,$ pointwise in $\Lambda_{1}(\alpha)$.

\end{itemize}
 \end{cor}
Since $H(\mathbb{D})$ and $\Lambda_{1}(n)$ are isomorphic, and all properties of interest remain invariant under isomorphism, we can state the following:

\begin{cor} Let $g\in H(\mathbb{D})$. The Rhaly operator $R_{g}: H(\mathbb{D})\to H(\mathbb{D})$  is
\begin{itemize}
\item[a.] topologizable,
\item[b.] m-topologizable if $\displaystyle \sup_{p\in \mathbb{N}}\|g\|_{p}<+\infty$.
\end{itemize}
The Rhaly operator $R_{g}:H(\mathbb{D})\to H(\mathbb{D})$ is also
\begin{itemize}
    \item[i.] power bounded,
    \item[ii.] mean ergodic,
    \item[iii.] uniformly mean ergodic,
    \item[iv.] Ces\`{a}ro bounded and $\displaystyle \lim_{n\to \infty} \frac{R^{n}_{g}}{n}=0,$ pointwise in $\Lambda_{\infty}(\alpha)$.
\end{itemize}
provided that  $\displaystyle \sup_{p\in \mathbb{N}}\|g\|_{p}\leq 1$.
\end{cor}
%%%%%%%%%%%%%%%%%%%%%%%%%%%%%%%%%%%%%%%%%%%%%%%%%%%%%%%%%%%%%%%%%%%%%%%%%%%%%%%%%%%%%%%%%%%%%%%%%%%%%%%%%%%%%%%%%%%%%%%%%%%%%%%%

%%%%%%%%%%%%%%%%%%%%%%%%%%%%%%%%%%%%%%%%%%%%%%%%%%%%%%%%%%%%%%%%%%%%%%%%%%%%%%%%%%%%%%%%%%%%%%%%%%%%%%%%%%%%%%%%%%%%%%%%%%%%%%%%%%%%%%%%%

% ------------------------------------------------------------------------
\end{document}